\documentclass[12pt]{article}
\usepackage{amssymb,latexsym}

\newcommand{\ben}{\begin{enumerate}}
\newcommand{\een}{\end{enumerate}}
\newcommand{\ble}{\begin{lem}}
\newcommand{\ele}{\end{lem}}
\newcommand{\bth}{\begin{thm}}
\renewcommand{\eth}{\end{thm}}
\newcommand{\bpr}{\begin{prop}}
\newcommand{\epr}{\end{prop}}
\newcommand{\bco}{\begin{cor}}
\newcommand{\eco}{\end{cor}}
\newcommand{\bcon}{\begin{conj}}
\newcommand{\econ}{\end{conj}}
\newcommand{\bde}{\begin{defn}}
\newcommand{\ede}{\end{defn}}
\newcommand{\bex}{\begin{exa}}
\newcommand{\eex}{\end{exa}}
\newcommand{\barr}{\begin{array}}
\newcommand{\earr}{\end{array}}
\newcommand{\btab}{\begin{tabular}}
\newcommand{\etab}{\end{tabular}}
\newcommand{\beq}{\begin{equation}}
\newcommand{\eeq}{\end{equation}}
\newcommand{\bea}{\begin{eqnarray*}}
\newcommand{\eea}{\end{eqnarray*}}
\newcommand{\bce}{\begin{center}}
\newcommand{\ece}{\end{center}}
\newcommand{\bpi}{\begin{picture}}
\newcommand{\epi}{\end{picture}}
\newcommand{\bfi}{\begin{figure} \begin{center}}
\newcommand{\efi}{\end{center} \end{figure}}

\newcommand{\bsl}{\begin{slide}{}}
\newcommand{\esl}{\end{slide}}

\newcommand{\bib}{thebibliography}
\newcommand{\pf}{{\bf Proof}\hspace{7pt}}

\newcommand{\Qqed}{\qquad\rule{1ex}{1ex}\medskip}

\newcommand{\hso}[1]{\hspace{-1pt}}

\newcommand{\case}[4]{\left\{\barr{ll}#1&\mbox{#2}\\#3&\mbox{#4}\earr\right.}

\def\<{\langle}
\def\>{\rangle}

\newcommand{\om}{\omega}

\newcommand{\De}{\Delta}

\newcommand{\bbN}{{\mathbb N}}

\newtheorem{thm}{Theorem}[section]
\newtheorem{prop}[thm]{Proposition}
\newtheorem{cor}[thm]{Corollary}
\newtheorem{lem}[thm]{Lemma}
\newtheorem{conj}[thm]{Conjecture}
\newtheorem{exa}[thm]{Example}

\begin{document}
\pagestyle{plain}

\title{On divisibility  of Narayana numbers by primes
}
\author{
Mikl\'os B\'ona \\[-5pt]
\small  Department of Mathematics, University of Florida\\[-5pt]
\small  Gainesville, FL 32611, USA, bona@math.ufl.edu
\\[5pt]
and\\[5pt]
Bruce E. Sagan\\[-5pt]
\small Department of Mathematics, Michigan State University\\[-5pt]
\small East Lansing, MI 48824-1027, USA, sagan@math.msu.edu
}

\date{\today}
\maketitle

\begin{abstract}
Using Kummer's Theorem, we give a necessary and sufficient condition
for a Narayana number to be divisible by a given prime.  We use this
to derive certain properties of the Narayana triangle.
\end{abstract}

\section{The main theorem}

Let $\bbN$ denote the nonnegative integers and let $k,n\in\bbN$.  The
{\it Narayana numbers\/}~\cite[A001263]{slo:ole} can be defined as
$$
N(n,k)=\frac{1}{n}{n\choose k}{n\choose k+1}
$$
where $0\le k<n$.  The Narayana numbers (in fact, a $q$-analogue of
them) were first studied by  MacMahon~\cite[Article 495]{mac:ca} and
were later rediscovered by Narayana~\cite{nar:tfp}.  They are closely
related to the {\it Catalan numbers\/}~\cite[A000108]{slo:ole}
$$
C_n=\frac{1}{n+1}{2n\choose n}
$$
and in fact $\sum_k N(n,k)= C_n$.   The
Narayana numbers can be arranged in a triangular array with $N(n,k)$
in row $n$ and column $k$ so that the row sums are the Catalan numbers.
Like the numbers $C_n$, the numbers
$N(n,k)$ have many combinatorial interpretations; see, for
example, the article of Sulanke~\cite{sul:nd}.

The main result of this note is a characterization of when $N(n,k)$ is
divisible by a given prime $p$.  To state it, we need some notation.  Let
$\De_p(n)=(n_i)$ denote the sequence of digits of $n$ in base $p$
so that $n=\sum_i n_i p^i$.  Similarly we define $\De_p(k)=(k_i)$.
If we are considering $k\le n$ then it will be convenient to 
extend the range of definition of $(k_i)$ so that both sequences have
the same length by setting $k_i=0$ if $p^i>k$.  The {\it order\/} of
$n$ modulo $p$ is the largest power of $p$ dividing $n$ and will be
denoted $\om_p(n)$.  As usual,
$k|n$ means that $k$ divides $n$.

Kummer's Theorem~\cite{kum:ear} gives a useful way of finding the
order of binomial coefficients.  For example, Knuth and
Wilf~\cite{kw:ppd} used it 
to find the highest power of a prime which divides a generalized
binomial coefficient.
\bth[Kummer]
Let $p$ be prime and let $\De_p(n)=(n_i)$, $\De_p(k)=(k_i)$. 
Then $\om_p {n\choose k}$ is the number of carries
in performing the addition $\De_p(k)+\De_p(n-k)$.  Equivalently, it is the
number of indices $i$ such that either $k_i>n_i$ or there exists an
index $j<i$ with $k_j>n_j$ and
$k_{j+1}=n_{j+1},\ldots,k_i=n_i$.\Qqed
\eth

Now everything is in place to state and prove our principal theorem.
\bth
\label{main}
Let $p$ be prime. 
Also let $\De_p(n)=(n_i)$, $\De_p(k)=(k_i)$ and $\om=\om_p(n)$. 
Then $p\nmid N(n,k)$ if and only if one of the two following conditions hold:
\ben 
\item When $p \nmid n$ we have
  \ben 
  \item $k_i\le n_i$ for all $i$, and
  \item $k_j < n_j$ where $j$ is the first index with $k_j\neq p-1$
  (if such an index exists).
  \een
\item When $p\mid n$ we have
  \ben 
  \item $k_i\le n_i$ for all $i>\om$, and
  \item $k_\om < n_\om$, and
  \item $k_0=k_1=\ldots=k_{\om-1}=
   \case{0}{if $p\mid k$;}{p-1}{if $p\nmid k$.}$
  \een
\een
\eth
\pf
First suppose that $p$ is not a divisor of $n$.  Then $p$ does not divide
$N(n,k)$ if and only if $p$ divides neither ${n\choose k}$ nor
${n\choose k+1}$.  By Kummer's Theorem this is equivalent to 
$k_i\le n_i$ and $(k+1)_i\le n_i$ for all $i$.  However, if $j$ is the
first index with $k_j\neq p-1$, then we have
$$
(k+1)_i=
{\left\{\barr{ll}0&\mbox{if $i<j$;}\\
(k)_i+1&\mbox{if $i=j$;}\\
(k)_i&\mbox{if $i>j$.}
\earr\right.}
$$
So these conditions can be distilled down to insisting that
$k_j<n_j$ in addition to $k_i\le n_i$ for all
other $i$.  

Now consider what happens when $p$ divides $n$.  Suppose first that
$p$ also divides $k$.  So $(n)_i=0$ for $i<\om$, which is a
nonempty set of indices, and $(k+1)_0=1$.  It follows
there are at least $\om$ carries in computing
$\De_p(k+1)+\De_p(n-k-1)$.  By Kummer's Theorem again, 
$\om_p{n\choose k+1}\ge\om$.  So $p$ does not divide $N(n,k)$ if
and only if it does not divide ${n\choose k}$ and $\om_p{n\choose k+1}=\om$.  
Applying Kummer's theorem  once more shows that this will happen exactly when
$k_i\le n_i$ for all $i$ with $k_\om<n_\om$.  So in particular $k_i=0$
for $i<\om$ since then $n_i=0$.  This completes the case when $p$
divides both $n$ and $k$.

Finally, suppose $p\mid n$ but $p\nmid k$.  Arguing as in the
previous paragraph, we see that $p$ is not a divisor of $N(n,k)$ if
and only if
$\om_p{n\choose k}=\om$ and $p$ does not divide ${n\choose k+1}$. 
But if $p$ is not a divisor of ${n\choose k+1}$ then, using Kummer's theorem,
 we must have $(k+1)_i=0$ for $i<\om$.  So $(k)_i=p-1$ for
$i<\om$.  Conditions 2(b) and (c) also follow as before.  This completes the
demonstration of the theorem.  \Qqed

\section{Applications}

It is well known that $C_n$ is odd if and only if $n=2^m-1$ for some
$m$.  For a combinatorial proof of this  which in fact establishes
$\om_2(C_n)$, see the article of Deutsch and Sagan~\cite{ds:ccm}.
Analogously, all the entries of the $n$th row of the Narayana triangle
are odd.  This is a special case of the following result.
\bco
\label{notdiv}
Let $p$ be prime and let $n=p^m-1$ for some $m\in\bbN$.  Then for 
all $k$, $0\le k\leq n-1$, we have $p\nmid N(n,k)$.

\eco
\pf
By Theorem~\ref{main} we just need to verify that 1(a) and (b) hold
for all $k$.  However, they must be true because
$n_i=p-1$ for all $i$.
\Qqed

We clearly can not have a row of the Narayana triangle where every
element is divisible by $p$ since $N(n,0)=N(n,n-1)=1$ for all $n$.
But we can ensure that every entry except the first and last is a
multiple of $p$.
\bco
\label{div}
Let $p$ be prime and let $n=p^m$ for some $m\in\bbN$.  Then $p\mid N(n,k)$
for $1\le k\le n-2$.
\eco
\pf
Suppose that $n=p^m$ and that $p$ does not divide $N(n,k)$.  If $p$
divides $k$, then 
condition 2(c) forces $k=0$.  If $p$ does not divide $k$, then the
same condition forces $k=n-1$.  So these are the only two numbers not
divisible by $p$ in the $n$th row of Narayana's triangle.
\Qqed

\section{Comments and Questions}

I.\  Clearly one could use the same techniques presented here to
determine $\om_p(N(n,k))$.  However, the cases become complicated enough
that it is unclear whether this would be an interesting thing to do.

\medskip

\noindent II.\  The characterization in Theorem~\ref{main} is involved enough
that it may be hopeless to ask for a combinatorial proof.  However,
there should be a combinatorial way to derive the simpler statements
in Corollaries~\ref{notdiv} and~\ref{div}, although we have not been
able to do so.  

As has already been mentioned, the order $\om_2(C_n)$
can be established by combinatorial means, specifically through the
use of group actions.  Unfortunately, the action used by
Deutsch and Sagan~\cite{ds:ccm} is not 
sufficiently refined to preserve the objects counted by $N(n,k)$.
For more information about how such methods can
be used to prove congruences, the reader can consult Sagan's
article~\cite{sag:cva} which also contains a survey of the literature.

Deutsch~\cite{deu:idp},
E\u{g}ecio\u{g}lu~\cite{ege:pcn}, and Simion and Ullman~\cite{su:sln}
have all found combinatorial ways to explain the fact that $C_n$ is
odd if and only if $n=2^m-1$ for some $m$.  Perhaps one or more of the
viewpoints in these papers could be adapted to the Narayana numbers.

\medskip

{\it Acknowledgment}.  We would like to thank Robert Sulanke for
valuable references and discussions, and to Neil White for bringing the
problem to our attention.

\begin{\bib}{99}

\bibitem{deu:idp}  E. Deutsch, An involution on Dyck paths and its
consequences, {\it Discrete Math.\/} {\bf 204} (1999), 163--166.

\bibitem{ds:ccm} 
\begin{flushleft}
E. Deutsch and B. E. Sagan, Congruences for Catalan
and Motzkin numbers and related sequences, preprint, available at
{\texttt http://www.math.msu.edu/$\sim$sagan/}.
\end{flushleft}

\bibitem{ege:pcn} \"O. E\u{g}ecio\u{g}lu, The parity of the Catalan
numbers via lattice paths, {\it Fibonacci Quart.\/}  {\bf 21} (1983)
65--66. 

\bibitem{kw:ppd} D. E. Knuth and H. S. Wilf, The power of a prime that
divices a generalized binomial coefficient, {\it J. Reine Angew.\
Math.\/} {\bf 396} (1989), 212--219.

\bibitem{kum:ear} E. E. Kummer, \"Uber die Erg\"anzungss\"atze zu den
allgemeinen Reciprocit\"ats\-gesetzen, {\it J. Reine Angew.\ Math.\/}
{\bf 44} (1852) 93--146.

\bibitem{mac:ca} P. A. MacMahon, {\it Combinatorial Analysis,
Vols. 1 and 2}, Cambridge University Press, 1915,1916; reprinted by
Chelsea, 1960.

\bibitem{nar:tfp}  T. V. Narayana, Sur les treillis form\'es par les
partitions d'une unties et leurs applications \`a la th\'eorie des
probabilit\'es, {\it C. R. Acad.\ Sci.\ Paris\/} {\bf 240} (1955),
1188--1189. 

\bibitem{sag:cva} B. E. Sagan, Congruences via Abelian groups,
{\it J. Number Theory} {\bf 20} (1985), 210--237.

\bibitem{su:sln} R. Simion and D. Ullman,  On the structure of the
lattice of noncrossing partitions, {\it Discrete Math.\/} {\bf 98}
(1991), 193--206.

\bibitem{slo:ole} 
\begin{flushleft}
N. J. A. Sloane, ``The On-Line Encyclopedia of
Integer Sequences,'' available at
{\texttt http://www.research.att.com/$\sim$njas/sequences/}.
\end{flushleft}

\bibitem{sul:nd} R. A. Sulanke, The Narayan distribution,  Special issue
on lattice path combinatorics and applications (Vienna, 1998).  {\it
J. Statist.\ Plann.\ Inference\/} {\bf 101} (2002), 311-326.

\end{\bib}

\btab{l}
\\
\hline
\\
{\it 2000 Mathematics Subject Classification:\/}  Primary 11B50;
Secondary 05A10, 11A07, 11B64.\\
{\it Keywords:\/} divisibility, Narayana numbers.
\\
\\
\hline
\\
(Concerned with sequences A001263 and A000108.)
\\
\\
\hline
\etab

\end{document}